\documentclass[english,12pt]{article}
\usepackage{amsfonts}
\usepackage{amsmath}
\usepackage[usenames]{color}
\usepackage[english]{babel}
\RequirePackage[cp1251]{inputenc}
\RequirePackage{srcltx}

\newtheorem{lemma}{Lemma}

\newtheorem{Theorem}{Theorem}

\DeclareMathOperator*{\supvrai}{sup \,\, vrai}
\begin{document}

\begin{center}{\bf ULYANOV INEQUALITIES FOR THE MIXED   MODULI OF SMOOTHNESS IN MIXED METRICS}
\footnote{ This research was partially funded by the Science Committee of the Ministry of Science and Higher Education of the Republic of Kazakhstan (Grant No. AP09260223).}
\end{center}

\vskip 0.2cm
\centerline{ B.V. Simonov$^1$, A.A. Jumabayeva$^2$ }
  \vskip 0.2cm

\centerline{\it $^{1}$  Volgograd State Technical University}
\centerline{\it  Volgograd, 400131, Russia}
\centerline{\it E-mail: simonov-b2002@yandex.ru}

\centerline{\it $^{2}$ L.N.Gumilyov Eurasian National University}
\centerline{\it  Satpayev 2, 010000 Astana, Kazakhstan}
\centerline{\it E-mail: ainurjumabay@gmail.com}

$\mathbb{}$ \ \ \ \ \  { \footnotesize

 $\mathbb{}$ \ \ \ \ \ \ \ In this paper we prove
Ulyanov-type inequalities between mixed moduli of smoothness of positive orders in different metrics. 
\vskip 0.1cm

 $\mathbb{}$ \ \ \ \ \ \ { \em Keywords}:
 Mixed  moduli of smoothness of positive order, Ulyanov type inequality. }

 \vskip 0.2cm

{ \bf  1. Introduction  \ }

\vskip 0.2cm

Let $L_{p_1p_2}, 1\leq p_i\le\infty, i=1,2 $
 be the set of measurable functions of two variables  $f(x_1, x_2), $ $ 2\pi$ -
 periodic in each variable, for which
$ \|f\|_{p_1 p_2}= || \{ ||f||_{p_1} \} ||_{p_2} < \infty, $  where
$$\|f\|_{p_i}=\left(\int\limits_{0}^{2\pi}
|f|^{p_i}dx_i\right)^{\frac1{p_i}}, \; \textnormal{if} \; 1 \leq
p_i < \infty, $$
$$\|f\|_{p_i}= \supvrai\limits_{0\leq x_i \leq 2 \pi}
 |f|, \; \textnormal{if} \;  p_i = \infty.$$

Let $L_{p_1p_2}^0$  be the space of functions $f\in L_{p_1p_2}$
 such that
 $\int\limits^{2\pi}_{0}f(x_1, x_2)dx_1 =0  \;    \textnormal{for almost all}  \; $ $  x_2$
 and
 $\int\limits^{2\pi}_{0}f(x_1, x_2)dx_2=0 \;   \textnormal{for almost all} \;   x_1.$

For the function $f\in L_{p_1p_2}$, we define the fractional  differences of positive order $\alpha_1$ and $\alpha_2$ with  steps $h_1$ and $h_2$ respectively, by variables $x_1$ and $x_2$ as follows:
$$
\Delta^{\alpha_1}_{h_1}(f)=\sum\limits^{\infty}_{\nu_1=0}(-1)^{\nu_1}
\left(^{\alpha_1}_{\nu_1}\right)f(x_1+(\alpha_1-\nu_1)h_1,x_2), $$
$$ \Delta^{\alpha_2}_{h_2}(f)=\sum\limits^{\infty}_{\nu_2=0}
(-1)^{\nu_2}\left(^{\alpha_2}_{\nu_2}\right)f(x_1,
x_2+(\alpha_2-\nu_2)h_2), $$
 where $\left(^{\alpha}_{\nu}\right)=1$
for  $\nu=0, \; \left(^{\alpha}_{\nu}\right)=\alpha$  for  $\nu=1
$,
$\left(^{\alpha}_{\nu}\right)=\frac{\alpha(\alpha-1)...(\alpha-\nu+1)}{\nu!}$
for  $\nu\geq 2. $

Denote (cf. \cite{Potapov2}) by $   \omega_{\alpha_1, \alpha_2}(f, \delta_1,
\delta_2)_{p_1p_2}$   the mixed modulus  of smoothness of positive orders $\alpha_1$ and $\alpha_2$, respectively, in the variables
$x_1$ and $x_2$ of a function  $f \in L_{p_1p_2}$,  that is,
$$   \omega_{\alpha_1, \alpha_2}(f, \delta_1, \delta_2 )_{p_1 p_2} =
\sup\limits_{|h_i| \leq \delta_i, i=1,2} ||
\Delta_{h_1}^{\alpha_1} \big(\Delta_{h_2}^{\alpha_2}(f) \big)
||_{p_1 p_2}.$$

Let $Y_{m_1, m_2}(f)_{p_1p_2}$ be the best approximation by the angle of the function $f \in L_{p_1p_2}$, i.e.
$$  Y_{m_1, m_2 }(f)_{p_1 p_2} =
\inf\limits_{T_{m_1, \infty }, T_{\infty, m_2 }} ||  f - T_{\infty,
m_2 } - T_{ m_1, \infty } ||_{p_1 p_2}, $$
where $ T_{m_1,\infty}(x_1, x_2)$ is a trigonometric polynomial of order at most
$m_1 (m_1 \in\mathbb{N} \cup \{0 \})$ with respect to $x_1$ and such that $T_{m_1, \infty} \in L_{p_1,p_2},$
 $ T_{\infty, m_2}(x_1, x_2)$ is a trigonometric polynomial of order at most $m_2 (m_2 \in \mathbb{N} \cup \{0 \})$ with respect to $x_2$ and such that $ T_{\infty, m_2} \in L_{p_1,p_2}$.

Let $  f^{(\rho_1, \rho_2)}(x_1,x_2) (\rho_1 \geq 0, \rho_2 \geq 0)$   be a derivative in the sense of Weyl of the function $f(x_1, x_2) \in L_{p_1,p_2}^0$ of order $\rho_1 $   with respect to $x_1$  and of order  $\rho_2 $ with respect to $x_2$ (see \cite{Besov1}, p. 238).

For a given number $q$, we will use the notation $q^*$, which will be defined as follows:
$q^*= q$ if $ 1 \leq q < \infty$ and $q^*=1$ if $q = \infty.$

For non-negative functionals $F(f,\delta_1,\delta_2)$ and $G(f,\delta_1,\delta_2)$ we write that $F(f,\delta_1,\delta_2)\ll G(f,\delta_1 ,\delta_2),$ if there exists a positive constant $C,$ independent of $f,\delta_1$ and $\delta_2$ such that $F(f,\delta_1,\delta_2)\leq C G(f,\delta_1 ,\delta_2).$ If $F(f,\delta_1,\delta_2)\ll G(f,\delta_1,\delta_2) $ and
$G(f,\delta_1,\delta_2)\ll F(f,\delta_1,\delta_2),$ then we will write that $F(f,\delta_1,\delta_2) \asymp G(f,\delta_1, \delta_2) $.

Denote by  $ \omega_{\alpha}(f,t)_p^{(1)} $  the moduli of smoothness of a fractional  order $\alpha, \; \alpha>0, $ of the function  $f\in L_p^{(1)}$, i.e.,
 $
 \omega_{\alpha}(f,t)^{(1)}_p =\sup\limits_{|h|\leq t}\|\Delta^{\alpha}_{h}f(x)\|^{(1)}_p,
$  where
 $$\Delta^{\alpha}_{h}f(x)=\sum\limits^{\infty}_{\nu=0}(-1)^{\nu}\left(^{\alpha}_{\nu}\right)f(x+(\alpha-\nu)h).$$

The  following $(p,q)$-inequality between moduli of smoothness in different metrics, nowadays
called sharp Ulyanov type inequalities, is known (see  \cite{ST3}):
$$
\omega_1(f,\delta)_q^{(1)} \ll \bigg( \int\limits_0^{\delta} \Big( t^{- \theta} \omega_1(f,t)_p^{(1)} \Big)^{q^*} \frac{dt}{t} \bigg)^{\frac{1}{q^*}},
 $$
where  $ 1 \leq p < q \leq \infty, \theta = \frac{1}{p} - \frac{1}{q}$.  For the
 Lebesgue spaces,
 Ulyanov type inequalities in the one-dimensional case have been studied for the moduli of smoothness of any positive order by many  authors, in particular, Ulyanov \cite{Ul1}, V.I. Kolyada
\cite{Kolyada}, W. Trebels \cite{Trebels}, O. Domingues \cite{DomTik},      B. Simonov \cite{ST3},   S. Tikhonov \cite{Tikhonov1}.

In the multidimensional case    Ulyanov type  inequalities were  considered by many authors. Inequalities  for the total moduli of smoothness of two variables were  obtained in \cite{DomTik}, \cite{Goldm}, \cite{KolTik1}, \cite{KolTik3},\cite{ST3}, and \cite{Tikhonov1}.
 Ulyanov type inequalities were   also proved for mixed moduli of smoothness in (\cite{Potapov2}-\cite{Potapov6}). The $(p,q)$ inequalities for the moduli of smoothness of derivatives in terms of moduli of smoothness of the function itself were  established in     \cite{Ainur}, \cite{KolTik1},
\cite{TikTre}.
Moreover, such inequalities were handled in various spaces in
\cite{DitTik}, \cite{GlazTik},  and \cite{GogOp}.

In this paper, we study Ulyanov-type inequalities for mixed moduli of smoothness.
For the case when $p_1=p_2=p$, the following relation between the mixed moduli of smoothness is known (\cite{Potapov2}-\cite{Potapov6}).
Let $ f \in L_{pp}^0, \alpha_i > 0, \delta_i \in (0,1) (i=1,2).$ Then

a) for $ 1 < p < q < \infty $
$$ \omega_{\alpha_1, \alpha_2}(f, \delta_1, \delta_2)_{qq}
\ll   \Bigg( \int\limits_0^{\delta_2}  \int\limits_0^{\delta_1}
\Big( t_1^{- \frac{1}{p} + \frac{1}{q}} t_2^{- \frac{1}{p} +
\frac{1}{q}} \omega_{ \alpha_1 + \frac{1}{p} - \frac{1}{q},
\alpha_2 + \frac{1}{p} - \frac{1}{q} } (f, t_1, t_2)_{pp} \Big)^{q}
\frac{d t_1}{t_1} \frac{d t_2}{t_2} \Bigg)^{\frac{1}{q}},
$$

b) for  $  1=p < q = \infty $
$$  \omega_{\alpha_1, \alpha_2} \big( f, \delta_1, \delta_2
\big)_{\infty \infty } \ll    \int\limits_{0}^{\delta_2}
\int\limits_{0}^{\delta_1}  (t_1 t_2)^{ - 1 } \omega_{\alpha_1 +
1, \alpha_2 + 1} (f, t_1, t_2)_{11} \frac{dt_1}{t_1}
\frac{dt_2}{t_2},$$\\

c) for $ 1 = p < q < \infty$
$$  \omega_{\alpha_1, \alpha_2} \big( f, \delta_1, \delta_2
\big)_{qq} \ll $$ $$  \ll \Big( \int\limits_{0}^{\delta_2 ( \log_2 \frac{2}{\delta_2})^{\frac{1}{ \alpha_2 q}} }  \int\limits_{0}^{\delta_1 (\log_2 \frac{2}{\delta_1})^{\frac{1}{ \alpha_1 q}} } \big( (t_1 t_2)^{ - 1 + \frac{1}{q} } \omega_{\alpha_1 + 1 - \frac{1}{q}, \alpha_2 + 1 - \frac{1}{q}} (f, t_1,
t_2)_{11} \big)^q  \frac{dt_1}{t_1} \frac{dt_2}{t_2} \Big)^{\frac{1}{q}},  $$

d) for $ 1 < p < q = \infty$
$$  \omega_{\alpha_1, \alpha_2} \big( f, \delta_1, \delta_2
\big)_{\infty \infty} \ll  \int\limits_{0}^{\delta_2 ( \log_2 \frac{2}{\delta_2})^{\frac{1}{ \alpha_2 p'}} }  \int\limits_{0}^{\delta_1 (\log_2 \frac{2}{\delta_1})^{\frac{1}{ \alpha_1 p'}} } (t_1 t_2)^{ - \frac{1}{p} }
\omega_{\alpha_1 + \frac{1}{p}, \alpha_2 + \frac{1}{p}} (f, t_1,
t_2)_{pp}  \frac{dt_1}{t_1} \frac{dt_2}{t_2}.  $$

Our main goal is to prove Ulyanov-type inequalities between mixed moduli of smoothness of positive orders, in particular,  we strengthen the above inequalities.
In this work, statements a) and b) are extended to the relation between mixed moduli of smoothness in the case when $p_1 \neq p_2.$
Let us formulate the main results.
\begin{Theorem}\label{th1}
 Let  $f \in L_{p_1 p_2}^0,$ where
 $ 1 < p_1 < q_1 < \infty  $ or   $ 1 = p_1 < q_1 = \infty$ and $ 1 < p_2 < q_2 < \infty $ or  $ 1 = p_2 < q_2 = \infty. $
Let for  $\alpha_i > 0, \delta_i \in (0,1), i = 1,2,$ we have
\begin{align}\label{eq1}
 &\omega_{\alpha_1, \alpha_2} \big( f, \delta_1,\delta_2 \big)_{q_1 q_2} \ll  \nonumber\\
&\Big( \int\limits_{0}^{\delta_2} \Big(
\int\limits_{0}^{\delta_1} \Big( t_1^{ - \frac{1}{p_1} + \frac{1}{q_1}} t_2^{ - \frac{1}{p_2} + \frac{1}{q_2} } \omega_{ \alpha_1 +  \frac{1}{p_1} - \frac{1}{q_1},\alpha_2  + \frac{1}{p_2} - \frac{1}{q_2}} (f, t_1, t_2)_{p_1 p_2} \Big)^{q_1^*} \frac{dt_1}{t_1} \Big)^{\frac{q_2^*}{q_1^*}} \frac{dt_2}{t_2} \Big)^{\frac{1}{q_2^*}}.
\end{align}
\end{Theorem}
{\bf Remark 1} {\it
Theorem \ref{th1} is optimal in the sense that there exists a function $f_0(x_1,x_2)$ such that for it in relation \eqref{eq1} the sign $\ll$ can be replaced by the sign $\asymp.$
}\\
{\bf Remark 2} {\it
In relation \eqref{eq1}, generally speaking, it is impossible to replace the sign $\ll$ with the sign $\asymp.$}

Let us now estimate the mixed moduli of smoothness of the derivative of functions in terms of the mixed moduli of smoothness of the function itself.

\begin{Theorem}\label{th2} Let  $f \in L_{p_1 p_2}^0  $, where
 $ 1 < p_1 < q_1 < \infty  $ or  $ 1 = p_1 < q_1 = \infty$ and $ 1 < p_2 < q_2 < \infty $ or $ 1 = p_2 < q_2 = \infty  $. Let for
$ \alpha_i > 0, \rho_i \geq 0, \delta_i \in (0,1), i = 1,2.$ we have
$$  \omega_{\alpha_1, \alpha_2} \big( f^{(\rho_1, \rho_2)}, \delta_1,\delta_2 \big)_{q_1 q_2} \ll $$
$$  \Big( \int\limits_{0}^{\delta_2} \Big(\int\limits_{0}^{\delta_1} \Big( t_1^{ - \rho_1 - \frac{1}{p_1} + \frac{1}{q_1}} t_2^{ - \rho_2 - \frac{1}{p_2} + \frac{1}{q_2} } \omega_{ \alpha_1 + \rho_1+ \frac{1}{p_1} - \frac{1}{q_1},\alpha_2 + \rho_2  + \frac{1}{p_2} - \frac{1}{q_2}} (f, t_1, t_2)_{p_1 p_2} \Big)^{q_1^*} \frac{dt_1}{t_1} \Big)^{\frac{q_2^*}{q_1^*}} \frac{dt_2}{t_2} \Big)^{\frac{1}{q_2^*}}. $$
\end{Theorem}

The  paper is organized as follows: in the next section, we provide some basic notions and results  that will be needed throughout the paper. In Sect. 3, we prove Theorems 1 and 2.

\vskip 0.2cm

 { \bf  2. Some auxiliary results }

 \vskip 0.2cm

 In this section, we give some notations and lemmas  which will be used in the proof of our main results.
Let $L_{p}^{(1)}, 1 \leq p \le \infty$ be the space of    $ 2\pi$  -
 periodic measurable functions for which
$ \|f\|_{p}^{(1)} < \infty, $ where
 $\|f\|_{p}^{(1)}=\left(\int\limits_{0}^{2\pi}
|f|^{p}dx\right)^{\frac{1}{p}}, \; \textnormal{when} \; 1 \leq
p < \infty, $ and  $\|f\|_{p}^{(1)}= \supvrai\limits_{0\leq x \leq 2 \pi}
 |f|, \; \textnormal{if} \;  p = \infty.$
Let
$L_{p}^{0(1)}$  be set of functions $f\in L_{p}^{(1)}$ such that   $\int\limits^{2\pi}_{0}f(x)dx=0.$

Let $ \sigma(f)^{(1)} $ be the Fourier series of the functions $f(x)$ i.e.
$$ \sigma(f)^{(1)} \equiv \sigma(f,x)^{(1)} \equiv \frac{a_0}{2} +
\sum\limits_{k=1}^\infty ( a_{k} \cos k x + b_{k} \sin k x) \equiv
\sum\limits_{k=0}^\infty A_{k}(x),$$
where $a_k$ and $b_k$ are the Fourier coefficients of the function $f(x)$.

Let $ \widetilde{f}(x) $ be  the conjugate function of $f(x)$ and
$f^{(\alpha)}(x)$ be the Weyl derivative of the function $f(x)$ of order $\alpha >0$.

$ T_n(\{\lambda_n\}; x) $ is the trigonometric polynomial $ T_n(x) = \sum\limits_{\nu=0}^n (c_\nu \cos \nu x + d_\nu \sin \nu x) $ transformed by the sequence $\{\lambda_n\}$  (\cite{JumSim}, \cite{Zigmund}) i.e.

$$
 T_n(\{\lambda_n\}; x) = \sum\limits_{\nu=0}^n \lambda_\nu (c_\nu \cos \nu x + d_\nu \sin \nu x). $$

Let
$V_{(1)m}(f) $  denote the de la Vall\'{e}e-Poussin sum   of the Fourier series of a function $f(x),$
$$ V_{ (1)m}(f)=\frac1{\pi}\int\limits^{2\pi}_{0}f(x+t)V_{m}^{2m}(t)dt, \ m = 0,1,2,.... $$

The two-dimensional de la Vall\'{e}e-Poussin sums of the Fourier series of a function  $f(x_1,x_2)$ are defined as follows
$$ V_{m_1, \infty}(f)=\frac1{\pi}\int\limits^{2\pi}_{0}f(x_1+t_1,x_2)V_{m_1}^{2m_1}(t_1)dt_1,$$
$$ V_{\infty, m_2}(f)=\frac1{\pi}\int\limits^{2\pi}_{0}f(x_1,x_2+t_2)V_{m_2}^{2m_2}(t_2)dt_2,$$
 $$ V_{m_1, m_2}(f)=\frac1{\pi^2}\int\limits^{2\pi}_{0}\int\limits^{2\pi}_{0}
 f(x_1+t_1,x_2+t_2)V_{m_1}^{2m_1}(t_1)V_{m_2}^{2m_2}(t_2)dt_1dt_2
 \ (m_i=0,1,2,..., \ i=1,2),$$
  where $V_0^0(t)=D_0(t),$ $ V_{n}^{2n}(t) = \frac{D_n(t) + ... + D_{2n-1}(t)}{n}, n=1,2,...,
\quad D_m(t)=\frac{\sin{(m+\frac12)t}}{2\sin{\frac{t}{2}}}, \;
m=0,1,2.... $

We need the following lemmas.
\begin{lemma}\label{lemma1}\cite{Potapov7}. Let $f \in L_{p_1p_2}^0, 1 \leq p_i \leq \infty,
\alpha_i > 0, n_i \in \mathbb{N}, i=1,2.$ Then

$\omega_{\alpha_1,\alpha_2}\Big(f,\frac{1}{n_1},\frac{1}{n_2}\Big)_{p_1p_2}
\asymp n^{- \alpha_1}_1n^{- \alpha_2}_2
\big\|V_{n_1,n_2}^{(\alpha_1,\alpha_2)}(f)\big\|_{p_1p_2} + n^{-
\alpha_1}_1 \big\|V_{n_1,\infty}^{(\alpha_1,0)}(f-
V_{\infty,n_2}(f))\big\|_{p_1p_2}$

$ + n^{- \alpha_2}_2
\big\|V_{\infty,n_2}^{(0,\alpha_2)}(f-V_{n_1,\infty}(f))\big\|_{p_1p_2}+
\big\|f-V_{n_1,\infty}(f)-V_{\infty,n_2}(f)+V_{n_1,n_2}(f)\big\|_{p_1p_2},$

$\|f-V_{n_1-1,\infty}(f)-V_{\infty,n_2-1}(f)+V_{n_1-1,n_2-1}(f)\|_{p_1p_2}
\ll Y_{n_1-1,n_2-1}(f)_{p_1p_2}.$

\end{lemma}

\begin{lemma}\label{lemma2}\cite{Potapov6}. Let $a_n \geq 0, b_n \geq 0,
\sum\limits_{k=1}^{n} a_k = a_n \gamma_n, 1 \leq p < \infty.$
Then
$$ \sum\limits_{k=1}^\infty \Big( \sum\limits_{n=k}^\infty b_n \Big)^{p}
\ll \sum\limits_{k=1}^\infty a_k ( b_k \gamma_k )^p. $$
\end{lemma}

 \begin{lemma}\label{lemma3}\cite{KolTik1,Potapov2,Potapov7}.
 Let $f \in L_{p_1p_2}^0$,$g\in L_{p_1p_2}^0$,  $1 \le p_i \le \infty,$  $\beta_i>\alpha_i> 0$,
  $n_i\in N$, $i=1,2.$ Then

 (1) $
\omega_{\alpha_1,\alpha_2}(f,\delta_1,0)_{p_1p_2}=\omega_{\alpha_1,\alpha_2}(f,0,\delta_2)_{p_1p_2}
=\omega_{\alpha_1,\alpha_2}(f,0,0)_{p_1p_2}=0.$

 (2) $
\omega_{\alpha_1,\alpha_2}(f+g,\delta_1,\delta_2)_{p_1p_2}\ll
     \omega_{\alpha_1,\alpha_2}(f,\delta_1,\delta_2)_{p_1p_2}+
     \omega_{\alpha_1,\alpha_2}(g,\delta_1,\delta_2)_{p_1p_2}.
$

  (3)
$ \omega_{\alpha_1,\alpha_2}(f,\delta_1,\delta_2)_{p_1p_2}\ll
\omega_{\alpha_1,\alpha_2}(f,t_1,t_2)_{p_1p_2}, $
 if  $ 0<\delta_i\le t_i, i=1,2$.

 (4) $
\frac{\omega_{\alpha_1,\alpha_2}(f,\delta_1,\delta_2)_{p_1p_2}}{\delta_1^{\alpha_1}\delta_2^{\alpha_2}}\ll
\frac{\omega_{\alpha_1,\alpha_2}(f,t_1,t_2)_{p_1p_2}}{t_1^{\alpha_1}t_2^{\alpha_2}},
$ if  $0<t_i\le\delta_i\le\pi, i=1,2$.

(5)
$\omega_{\alpha_1,\alpha_2}(f,\lambda_1\delta_1,\lambda_2\delta_2)_{p_1p_2}\ll\,
(\lambda_1+1)^{\alpha_1}(\lambda_2+1)^{\alpha_2}\,
\omega_{\alpha_1,\alpha_2}(f,\delta_1,\delta_2)_{p_1p_2},$
 if  $ \lambda_i>0,i=1,2.$

(6) $ Y_{n_1-1,n_2-1}(f)_{p_1p_2} \ll \omega_{\alpha_1,\alpha_2}
(f,\frac{1}{n_1},\frac{1}{n_2})_{p_1p_2} \ll $

$\frac{1}{n_1^{\alpha_1}}\frac{1}{n_2^{\alpha_2}}\sum\limits_{v_1=1}^{n_1+1}\sum\limits_{v_2=1}^{n_2+1}
v_1^{\alpha_1-1}v_2^{\alpha_2-1}Y_{v_1-1,v_2-1}(f)_{p_1p_2}. $

 (7)
$\omega_{\beta_1,\beta_2}(f,\delta_1,\delta_2)_{p_1p_2}\ll
\omega_{\alpha_1,\alpha_2}(f,\delta_1,\delta_2)_{p_1p_2}.$

 (8)
$\frac{\omega_{\alpha_1,\alpha_2}(f,\delta_1,\delta_2)_{p_1p_2}}
{\delta_1^{\alpha_1}\delta_2^{\alpha_2}}\ll
\frac{\omega_{\beta_1,\beta_2}(f,\delta_1,\delta_2)_{p_1p_2}}
{\delta_1^{\beta_1}\delta_2^{\beta_2}}$, if $0<\delta_i\le\pi$,
$i=1,2$.

(9) $ \delta_1^{-r_1} \omega_{\alpha_1 + r_1, \alpha_2}(f,
\delta_1, \delta_2)_{p_1 p_2} \ll \omega_{\alpha_1, \alpha_2}(f^{(r_1, 0)},
\delta_1, \delta_2)_{p_1 p_2} \ll $

$ \ll \int\limits_0^{\delta_1} t_1^{-r_1}
\omega_{\alpha_1 + r_1, \alpha_2}(f, t_1, \delta_2)_{p_1 p_2} \frac{d t_1}{t_1}, $
if $0< \delta_1 < \pi$ and $ f^{(r_1,0)} \in L_{p_1 p_2}^0, $

$ \delta_2^{-r_2} \omega_{\alpha_1, \alpha_2 + r_2}(f,  \delta_1,
\delta_2)_{p_1 p_2} \ll \omega_{\alpha_1, \alpha_2}(f^{(0, r_2)},  \delta_1,
\delta_2)_{p_1 p_2} \ll $

$ \ll \int\limits_0^{\delta_2} t_2^{-r_2}
\omega_{\alpha_1, \alpha_2 + r_2}(f, \delta_1, t_2)_{p_1 p_2} \frac{d t_2}{t_2}, $
if $0< \delta_2 < \pi$ and $f^{(0, r_2)} \in L_{p_1p_2}^0, $

$ \delta_2^{-r_2} \delta_1^{-r_1} \omega_{\alpha_1+r_1, \alpha_2 + r_2}(f, \delta_1,
\delta_2)_{p_1 p_2} \ll \omega_{\alpha_1, \alpha_2}(f^{(r_1, r_2)}, \delta_1,
\delta_2)_{p_1 p_2} \ll $

$ \ll \int\limits_0^{\delta_2}  \int\limits_0^{\delta_1} t_2^{-r_2} t_1^{-r_1}
\omega_{\alpha_1+r_1, \alpha_2 + r_2}(f, t_1, t_2)_{p_1 p_2} \frac{d t_2}{t_2}, $
if    $0< \delta_i < \pi$ $(i=1,2)  $ and $f^{(r_1, r_2)} \in L_{p_1p_2}^0. $
\end{lemma}

 \begin{lemma}\label{lemma4}\cite{Potapov8}. Let $ f \in L_{p_1 p_2}^0, 1 \leq p_i < q_i
\leq \infty, n_i = 0,1,2,..., i =1,2.$ Then
$$ Y_{2^{n_1}-1, 2^{n_2}-1}(f)_{q_1 q_2} \ll \Bigg( \sum\limits_{\nu_2 =n_2}^\infty \Big(
\sum\limits_{\nu_1 =n_1}^\infty \Big( 2^{\nu_1  \big( \frac{1}{p_1} - \frac{1}{q_1} \big)} 2^{ \nu_2 \big( \frac{1}{p_2} - \frac{1}{q_2} \big)  } Y_{2^{\nu_1} -1, 2^{\nu_2 }-1}(f)_{p_1 p_2} \Big)^{q_1^*}  \Big)^{\frac{q_2^*}{q_1^*}} \Bigg)^{\frac{1}{q_2^*}}.
$$
\end{lemma}

 \begin{lemma}\label{lemma5}\cite{PS8}. (a) Let the function $T_{ \infty, 2^{n_2} } \in L_{p_1 p_2},$ $
1 \leq p_i  \leq \infty (i=1,2) $ be a trigonometric polynomial of order $2^{n_2}$ in variable $x_2,$ $ 1 \leq p_2 < q_2
\leq \infty.$ Then
$$ \Big\| \sum\limits_{n_2=M_1}^{M_2} T_{ \infty, 2^{n_2}} \Big\|_{p_1q_2} \leq C \left(
\sum\limits_{n_2=M_1}^{M_2}  2^{n_2 q_2^*(\frac{1}{p_2} -
\frac{1}{q_2})} \| T_{ \infty, 2^{n_2}} \|_{p_1p_2}^{q_2^*}
\right)^{1/q_2^*},
$$
where constant C does not depend on $N_1, N_2$ and $T_{2^{n_1}, \infty}.$

(b) Let the function $T_{2^{n_1}, \infty } \in L_{p_1 p_2},$ $
1 \leq p_i  \leq \infty (i=1,2) $ be a trigonometric polynomial of order $2^{n_1}$ in variable $x_1,$ $ 1 \leq p_1 < q_1
\leq \infty. $  Then
$$ \Big\| \sum\limits_{n_1=N_1}^{N_2} T_{2^{n_1}, \infty} \Big\|_{q_1p_2} \leq C \left(
\sum\limits_{n_1=N_1}^{N_2}  2^{n_1 q_1^*(\frac{1}{p_1} -
\frac{1}{q_1})} \| T_{2^{n_1}, \infty} \|_{p_1p_2}^{q_1^*}
\right)^{1/q_1^*},
$$
where constant C does not depend on  $N_1, N_2$ and $T_{2^{n_1}, \infty}.$
\end{lemma}

 \begin{lemma}\label{lemma6}\cite{PS6}.  Let $ f \in L_{p_1p_2}^0, 1 \leq p_i \leq
\infty, n_i =0,1,2, ...,  i=1,2. $ Then

(a) $  \| V_{\infty 2^{n_2}}^{(0, \alpha_2)} (f - V_{2^{n_1} \infty}(f)) \|_{p_1p_2} \ll
 \| \sum\limits_{\nu_1= n_1}^\infty V_{\infty 2^{n_2}}^{(0, \alpha_2)} ( V_{2^{\nu_1} \infty}(f) - V_{2^{\nu_1 +1}\infty}(f) ) \|_{p_1p_2}, $

(b) $   \| V_{ 2^{n_1} \infty}^{( \alpha_1, 0)}  (f - V_{\infty 2^{n_2} }(f)) \|_{p_1p_2} \ll
  \| \sum\limits_{\nu_2= n_2}^\infty V_{ 2^{n_1} \infty}^{( \alpha_1, 0)} ( V_{ \infty 2^{\nu_2}}(f) -V_{\infty 2^{\nu_2 +1}}(f) ) \|_{p_1p_2}. $

\end{lemma}

 \begin{lemma}\label{lemma7}  (\cite{Zigmund}, p. 213). Let $ f(x) \in L_p^{0(1)}, 1 < p < q < \infty,
\theta = \frac{1}{p}-\frac{1}{q}, \alpha \geq 0$. Then
$$\parallel f^{(\alpha)} \parallel_q^{(1)} \ll \parallel f^{(\alpha + \theta)} \parallel_p^{(1)}.$$
\end{lemma}

 \begin{lemma}\label{lemma8}\cite{PS6}.
(a) Let the function $T_{n_1, \infty}(x_1,x_2) \in L_{p_1 p_2}^{0}, 1 = p_1 < q_1=\infty,
1 \leq p_2 \leq \infty$ be a trigonometric polynomial of order at most $n_1 (n_1 \in \mathbb{N}),$ in the variable  $x_1.$ Then
$$ \|T_{n_1, \infty}\|_{\infty p_2} \leq C_1 \| T_{n_1, \infty}^{(1,0)} \|_{1, p_2}, $$
where the constant $C_1$ does not depend on  $T_{n_1, \infty}.$

(b) Let the function $T_{\infty, n_2}(x_1,x_2) \in L_{p_1 p_2}^{0}, 1 = p_2 < q_2=\infty,
1 \leq p_1 \leq \infty,$ be a trigonometric polynomial of order at most  $n_2 (n_2 \in \mathbb{N}),$ in the variable  $x_2.$ Then
$$ \|T_{\infty, n_2}\|_{p_1 \infty} \leq C_2 \| T_{\infty, n_2}^{(0,1)} \|_{p_1, 1}, $$
where the constant $C_2$ does not depend on  $T_{\infty, n_2}.$
\end{lemma}

 \begin{lemma}\label{lemma9}(\cite{PST4}, p. 13). Let  $\alpha >0, f(x) \in L_1^{0(1)}, \sigma(f)^{(1)} = \sum\limits_{k=1}^\infty A_k(x),$ where \\ $A_k(x) = a_k \cos kx + b_k \sin kx.$ Let the function $\varphi(x)$ be such that $ \sigma(\varphi)^{(1)} = \sum\limits_{k=1}^\infty k^\alpha A_k(x). $ Then

\centerline{$ \varphi(x) = \cos \frac{ \pi \alpha}{2} f^{(\alpha)} (x) + \sin \frac{\pi \alpha}{2} \widetilde{f}^{(\alpha)}(x). $}
\end{lemma}

 \begin{lemma}\label{lemma10}\cite{PS5}. Let $0 < \alpha < 1, n \in \mathbb{N}, 0 < |x| \leq \pi.$ Then

$ |\sum\limits_{\nu=1}^n \frac{ \cos \nu x}{ \nu^\alpha}| \leq C(\alpha) \cdot |x|^{ \alpha -1}, $

$ |\sum\limits_{\nu=1}^n \frac{ \sin \nu x}{ \nu^\alpha}| \leq C(\alpha) \cdot |x|^{ \alpha -1}, $
\\ where the positive constant $C(\alpha)$ does not depend on $n$ and $x,$ but depends only on $\alpha.$
\end{lemma}

 \begin{lemma}\label{lemma11} (\cite{Besov1}, p. 31). (Hardy-Littlewood inequality). Let $
1 < p < q < \infty, $  $ \mu = 1 - \frac{1}{p} + \frac{1}{q}, $ $f \in L_p^{(1)}(\mathbb{R}_1),$ $ \mathbb{R}_1= (- \infty, + \infty), $ that is $ \|f\|_{L_p^{(1)}(\mathbb{R}_1)} = \Big(  \int\limits_{- \infty}^\infty |f(x)|^p dx  \Big)^{\frac{1}{p}} < \infty,$
$$ J(x) = \int\limits_{\mathbb{R}_1} f(y) |y-x|^{- \mu} dy. $$
Then the inequality is true
$$ \|J\|_{L_q^{(1)} (\mathbb{R}_1)} \leq K(p,q) \|f\|_{{L_p^{(1)}} (\mathbb{R}_1)}, $$
where the constant $K(p,q)$ depends only on $p$ and $q.$
\end{lemma}

 \begin{lemma}\label{lemma12}  Let   $ 1 < p_1 < q_1 < \infty,
1 \leq p_2 \leq \infty.$ Let the function $T_{n_1, \infty}(x_1,x_2) \in L_{p_1 p_2}^{0}$ be a trigonometric polynomial of order at most $n_1 (n_1 \in \mathbb{N}),$ in $x_1.$ Then
$$ \|T_{n_1, \infty}\|_{q_1 p_2} \leq C_1 \| T_{n_1, \infty}^{(\frac{1}{p_1}-\frac{1}{q_1},0)} \|_{p_1, p_2}, $$
where the constant $C_1$ does not depend on $T_{n_1, \infty}.$
\end{lemma}

{ \bf   Proof. }
Let $1 \leq p_2 < \infty. $ Let us estimate
$$\|T_{n_1, \infty}\|_{q_1 p_2} =   \Big( \int\limits_{0}^{2 \pi} \Big( \Big( \int\limits_{0}^{2 \pi} \Big| T_{n_1, \infty}(x_1, x_2) \Big|^{q_1} dx_1 \Big)^{\frac{1}{q_1}}  \Big)^{p_2} dx_2 \Big)^{\frac{1}{p_2}}. $$
For almost all fixed $x_2^0$, the function $T_{n_1, \infty}(x_1, x_2^0) = P_{n_1}(x_1)$ is a function of one variable $x_1.$ Applying Lemma \ref{lemma7}, we get
$$\Big( \int\limits_{0}^{2 \pi} \Big| P_{n_1}(x_1) \Big|^{q_1} dx_1 \Big)^{\frac{1}{q_1}} \ll \Big( \int\limits_{0}^{2 \pi} \Big| P_{n_1}^{(\frac{1}{p_1} - \frac{1}{q_1})}(x_1) \Big|^{p_1} dx_1 \Big)^{\frac{1}{p_1}}.$$
Then the following  inequality is true
$$  \Big( \int\limits_{0}^{2 \pi} \Big( \Big( \int\limits_{0}^{2 \pi} \Big| T_{n_1, \infty}(x_1, x_2) \Big|^{q_1} dx_1 \Big)^{\frac{1}{q_1}}  \Big)^{p_2} dx_2 \Big)^{\frac{1}{p_2}}
$$
$$\ll \Big( \int\limits_{0}^{2 \pi} \Big( \Big( \int\limits_{0}^{2 \pi} \Big| T_{n_1, \infty}^{(\frac{1}{p_1} - \frac{1}{q_1},0)}(x_1, x_2) \Big|^{p_1} dx_1 \Big)^{\frac{1}{p_1}}  \Big)^{p_2} dx_2 \Big)^{\frac{1}{p_2}}, $$
whence follows Lemma \ref{lemma12}. Lemma \ref{lemma12}  is verified similarly in the case when $p_2 = \infty.$ Lemma \ref{lemma12} is proved.

 \begin{lemma}\label{lemma13} Let $ 1 < p_2 < q_2 < \infty,1 \leq p_1 \leq \infty$. Let the function  $T_{\infty, n_2}(x_1,x_2) \in L_{p_1 p_2}^{0}$
be a trigonometric polynomial of order at most $n_2 (n_2 \in \mathbb{N}),$ in $x_2.$ Then
$$ \|T_{\infty, n_2}\|_{p_1 q_2} \leq C_2 \| T_{\infty, n_2}^{(0, \frac{1}{p_2} - \frac{1}{q_2})} \|_{p_1, p_2}, $$
where the constant $C_2$ does not depend on  $T_{\infty, n_2}.$
\end{lemma}

{ \bf   Proof. }
For almost all fixed $x_1^0$, the function $P_{n_2}(x_2) = T_{\infty, n_2}(x_1^0, x_2)$ is a function of one variable $x_2.$ Then

 \centerline{$  P_{n_2}(x_2) = \sum\limits_{\nu=1}^{n_2} (c_\nu \cos \nu x_2 + d_\nu \sin \nu x_2). $}
Let  $ \theta_2 = \frac{1}{p_2} - \frac{1}{q_2}. $ Then, applying Lemma \ref{lemma9}, we have
\begin{align*}
&  P_{n_2}(x_2) = \sum\limits_{\nu=1}^{n_2} (c_\nu \cos \nu x_2 + d_\nu \sin \nu x_2) = \sum\limits_{\nu=1}^{n_2} \frac{1}{\nu^{\theta_2}} \nu^{\theta_2} (c_\nu \cos \nu x_2 + d_\nu \sin \nu x_2) \\
&= \sum\limits_{\nu=1}^{n_2} \frac{1}{\nu^{\theta_2}} \Big( \cos \frac{\theta_2 \pi}{2} (c_\nu \cos \nu x_2 + d_\nu \sin \nu x_2)^{(\theta_2)} + \sin \frac{\theta_2 \pi}{2}
\widetilde{(c_\nu \cos \nu x_2 + d_\nu \sin \nu x_2)}^{(\theta_2)}\Big) \\
&= \cos \frac{\theta_2\pi}{2}\big(P_{n_2}(x_2)(\{\frac{1}{m^{\theta_2}}\};x_2)\big)^{(\theta_2)} + \sin \frac{\theta_2\pi}{2}\big(\widetilde{P}_{n_2}(\{\frac{1}{m^{\theta_2}}\};x_2)\big)^{(\theta_2)}\\
&=  \cos \frac{\theta_2\pi}{2} \frac{1}{\pi} \int\limits_{-\pi}^{\pi} \Big(P_{n_2}(t_2)\Big)^{(\theta_2)}
\sum\limits_{\nu=1}^{n_2} \frac{1}{\nu^{\theta_2}} \cos \nu (t_2-x_2) dt_2  \\
&\;\;\;\;\;\;\;\;\;\;\;\;\;- \sin \frac{\theta_2\pi}{2} \frac{1}{\pi} \int\limits_{-\pi}^{\pi} \Big(P_{n_2}(t_2)\Big)^{(\theta_2)}
\sum\limits_{\nu=1}^{n_2} \frac{1}{\nu^{\theta_2}} \sin \nu (t_2-x_2) dt_2.
\end{align*}

Now, using Lemma \ref{lemma10}, for almost all $x_2$ we obtain
\begin{align*}
&  |P_{n_2}(x_2)| \leq   \frac{1}{\pi} \int\limits_{-\pi}^{\pi} \Big| \Big(P_{n_2}(t_2)\Big)^{(\theta_2)} \Big|
\Big|\sum\limits_{\nu=1}^{n_2+1} \frac{1}{\nu^{\theta_2}} \cos \nu (t_2-x_2) \Big| dt_2  \\
&\;\;\;\;\;\;\;\;\;\;  + \frac{1}{\pi} \int\limits_{-\pi}^{\pi} \Big| \Big(P_{n_2}(t_2)\Big)^{(\theta_2)} \Big|
\Big|\sum\limits_{\nu=1}^{n_2} \frac{1}{\nu^{\theta_2}} \sin \nu (t_2-x_2) \Big| dt_2   \\
&  \leq C(\theta_2)   \frac{1}{\pi} \int\limits_{-\pi}^{\pi} \Big| \Big(P_{n_2}(t_2)\Big)^{(\theta_2)} \Big|\Big| t_2-x_2  \Big|^{\theta_2 -1} dt_2 + C (\theta_2) \frac{1}{\pi} \int\limits_{-\pi}^{\pi} \Big| \Big(P_{n_2}(t_2)\Big)^{(\theta_2)} \Big|\Big|t_2-x_2 \Big|^{\theta_2-1} dt_2
\\
&  \leq C_1(\theta_2)    \int\limits_{-\pi}^{\pi} \Big| \Big(P_{n_2}(t_2)\Big)^{(\theta_2)} \Big|\Big| t_2-x_2  \Big|^{\theta_2 -1} dt_2.
\end{align*}

Thus, for almost all $x_2$ the following inequality holds
$$ |P_{n_2}(x_2)| \leq C(p_2,q_2) \int\limits_{-\pi}^{\pi} \Big| P_{n_2}^{(  \theta_2)}(t_2) \Big|\Big| t_2-x_2  \Big|^{\theta_2 -1} dt_2.$$

Therefore, for almost all $x_1$ and $x_2$, the inequality
\begin{equation}\label{eq2}
 |T_{\infty, n_2}(x_1, x_2)|  \leq C(p_2,q_2) \int\limits_{-\pi}^{\pi} \Big|  T_{\infty, n_2}^{(0,  \theta_2)}(x_1,t_2) \Big|\Big| t_2-x_2  \Big|^{\theta_2 -1} dt_2,
  \end{equation}
is true, where the constant $C(p_2,q_2)$ depends only on $p_2$ and $q_2.$

If $p_1 < \infty,$ then,  using \eqref{eq2},  we obtain
\begin{align*}
&  \|T_{\infty, n_2}  \|_{p_1 q_2}=   \Big( \int\limits_{-\pi}^\pi \Big(  \Big\{ \int\limits_{-\pi}^\pi \Big|T_{\infty, n_2}(x_1,x_2)  \Big|^{p_1} dx_1 \Big\}^{\frac{1}{p_1}} \Big)^{q_2} dx_2 \Big)^{\frac{1}{q_2}} \leq \\
& \leq C(p_2,q_2)    \Big( \int\limits_{-\pi}^\pi \Big(  \Big\{ \int\limits_{-\pi}^\pi \Big|\int\limits_{-\pi}^{\pi} \Big|  T_{\infty, n_2}^{(0,  \theta_2)}(x_1,t_2) \Big|\Big| t_2-x_2  \Big|^{\theta_2 -1} dt_2 \Big|^{p_1} dx_1 \Big\}^{\frac{1}{p_1}} \Big)^{q_2} dx_2 \Big)^{\frac{1}{q_2}}.
\end{align*}
Using the Minkowski inequality, we have
\begin{align*}
&  \|T_{\infty, n_2}  \|_{p_1 q_2}  \leq \\
& \leq C(p_2,q_2)  \Big( \int\limits_{-\pi}^\pi \Big(   \int\limits_{-\pi}^\pi \Big\{\int\limits_{-\pi}^{\pi} \Big|  T_{\infty, n_2}^{(0,  \theta_2)}(x_1, x_2)(x_1,t_2) \Big|^{p_1}dx_1 \Big\}^{\frac{1}{p_1}}  \Big| t_2-x_2  \Big|^{\theta_2 -1} dt_2  \Big)^{q_2} dx_2 \Big)^{\frac{1}{q_2}}.
\end{align*}

We introduce the function $F(t_2)$, which is equal to $\Big\{\int\limits_{-\pi}^{\pi} \Big|   T_{\infty, n_2}^{(0,  \theta_2)}(x_1,t_2) \Big|^{p_1}dx_1 \Big\}^{\frac{1}{p_1}},$
 if $t_2 \in [-\pi, \pi]$ and

 $$ \|T_{\infty, n_2}  \|_{p_1 q_2}  \leq C(p_2,q_2)  \Big( \int\limits_{\mathbb{R}_1} \Big(   \int\limits_{\mathbb{R}_1} F(t_2)  \Big| t_2-x_2  \Big|^{\theta_2 -1} dt_2  \Big)^{q_2} dx_2 \Big)^{\frac{1}{q_2}}.$$
 Therefore, from Lemma  11, we deduce that
 $$ \|T_{\infty, n_2}  \|_{p_1 q_2}   \leq C_1(p_2,q_2)  \|F\|_{L_{p_2}^{(1)}(\mathbb{R}_1)} = C_1(p_2,q_2)  \| T_{\infty, n_2}^{(0,  \theta_2)} \|_{p_1 p_2}.$$
If $p_1= \infty$, then, applying inequality \eqref{eq2}, we get
   \begin{align*}
  & \|T_{\infty, n_2}  \|_{p_1 q_2} =    \Big(  \int\limits_{- \pi}^\pi \Big[ \supvrai\limits_{0\leq x_1 \leq 2 \pi} | T_{\infty, n_2}(x_1,x_2) | \Big]^{q_2} dx_2  \Big)^{\frac{1}{q_2}} \\
&  \ll   \Big(  \int\limits_{- \pi}^\pi \Big[ \supvrai\limits_{0\leq x_1 \leq 2 \pi} | \int\limits_{-\pi}^{\pi} \Big|  T_{\infty, n_2}^{(0,  \theta_2)}(x_1,t_2) \Big|\Big| t_2-x_2  \Big|^{\theta_2 -1} dt_2 | \Big]^{q_2} dx_2  \Big)^{\frac{1}{q_2}} \\
&  \ll   \Big(  \int\limits_{- \pi}^\pi \Big[  \int\limits_{-\pi}^{\pi} \supvrai\limits_{0\leq x_1 \leq 2 \pi} | \Big|  T_{\infty, n_2}^{(0,  \theta_2)}(x_1,t_2) \Big|\Big| t_2-x_2  \Big|^{\theta_2 -1} dt_2 | \Big]^{q_2} dx_2  \Big)^{\frac{1}{q_2}}.
 \end{align*}
 Consider the function $F_1(t_2)$, which is equal to $ \supvrai\limits_{0\leq x_1 \leq 2 \pi} | T_{\infty, n_2}^{(0, \theta_2)}(x_1,t_2)|, $  if $t_2 \in [-\pi, \pi]$ and is zero otherwise. Then we obtain
 $$ \|T_{\infty, n_2}  \|_{p_1 q_2}  \leq C(p_2,q_2)   \Big( \int\limits_{\mathbb{R}_1} \Big(   \int\limits_{\mathbb{R}_1} F_1(t_2)  \Big| t_2-x_2  \Big|^{\theta_2 -1} dt_2  \Big)^{q_2} dx_2 \Big)^{\frac{1}{q_2}}.$$
 Applying Lemma \ref{lemma11}, we conclude that
  $$ \|T_{\infty, n_2}  \|_{p_1 q_2}  \leq \|F_1\|_{L_{p_2}^{(1)}(\mathbb{R}_1)} \leq C_1(p_2,q_2)  \| T_{\infty, n_2}^{(0,  \theta_2)}\|_{p_1 p_2}.$$
Thus, we have shown that
$$ \|T_{\infty, n_2}  \|_{p_1 q_2}  \leq  C_1(p_2,q_2)  \| T_{\infty, n_2}^{(0, \theta_2)}\|_{p_1 p_2}.$$
Lemma \ref{lemma13} is proved.

 \begin{lemma}\label{lemma14}(\cite{PST4},\cite{Tikhonov3}). (a) Let $ 1 < p < q < \infty, \ \delta \in (0,1).$
 Consider the function  $ g_1(x) = \sum\limits_{\nu =0}^\infty \frac{( \nu +1)^\beta}{2^{\nu \alpha}} \cos 2^\nu x, $   where $\beta > -\frac{1}{2}, \alpha > \frac{1}{p} - \frac{1}{q}.$  Then
 $$ \omega_\alpha(g_1, \delta)_q^{(1)} \asymp \delta^\alpha \bigg(
\log_2 \frac{2}{\delta} \bigg)^{\beta + \frac{1}{2}},$$
$$ \Bigg(
\int\limits_0^{\delta} \Big[ t^{ - \theta} \omega_{\alpha + \frac{1}{p} - \frac{1}{q}}(g_1,t)_p^{(1)} \Big]^{q}
\frac{dt}{t} \Bigg)^{\frac{1}{q}} \asymp \delta^{ \alpha - \frac{1}{p} + \frac{1}{q}
} \bigg( \log_2 \frac{2}{\delta} \bigg)^{\beta }.$$
(b) Let $ 1 = p < q = \infty, \ \delta \in (0,1).$
  Consider the function   $ g_2(x) = \sum\limits_{\nu =0}^\infty \frac{( \nu +1)^\beta}{2^{\nu \alpha}} \cos  \big( 2^\nu x - \frac{ \pi \alpha}{2} \big), $
  where $\beta > - 1, \alpha > 1.$  Then
 $$ \omega_\alpha(g_2, \delta)_\infty^{(1)} \asymp \delta^\alpha \bigg(
\log_2 \frac{2}{\delta} \bigg)^{\beta +1},$$
$$ \Bigg(
\int\limits_0^{\delta}  t^{ - 1} \omega_{\alpha + 1}(g_2,t)_1^{(1)}
\frac{dt}{t} \Bigg)^{\frac{1}{q}} \asymp \delta^{ \alpha - 1
} \bigg( \log_2 \frac{2}{\delta} \bigg)^{\beta }.$$
\end{lemma}

 \vskip 0.2cm

 { \bf  3. Proof of main results }

{ \bf  Proof of Theorems \ref{th1}. } For every $\delta_i \in (0,1)$, take  a non-negative integer $n_i$ such that  $
\frac{1}{2^{n_i +1}} \leq \delta_i < \frac{1}{2^{n_i}}, i =1,2. $
Using Lemma \ref{lemma3}, we have
$$ I = \omega_{\alpha_1, \alpha_2} \big( f, \delta_1,
\delta_2 \big)_{q_1 q_2} \ll \omega_{\alpha_1, \alpha_2} \Big( f,
\frac{1}{2^{n_1}}, \frac{1}{2^{n_2}} \Big)_{q_1 q_2}.  $$
Applying Lemma \ref{lemma1}, we derive the following
\begin{align*}
 I &\ll 2^{-n_1 \alpha_1 - n_2 \alpha_2} \| V_{2^{n_1},2^{n_2}}^{(\alpha_1, \alpha_2)}(f) \|_{q_1 q_2} + 2^{-n_1 \alpha_1 }\| V_{2^{n_1}, \infty}^{(\alpha_1, 0)}(f - V_{\infty, 2^{n_2}}(f)) \|_{q_1 q_2} \\
 &+ 2^{-n_2 \alpha_2 } \| V_{ \infty, 2^{n_2}}^{(0, \alpha_2)} (f- V_{2^{n_1}, \infty }(f) )  \|_{q_1 q_2} + \| f - V_{2^{n_1}, \infty}(f) -V_{ \infty, 2^{n_2}}(f) + V_{ 2^{n_1}, 2^{n_2}}(f) \|_{q_1 q_2}\\
& \equiv  2^{- n_1 \alpha_1 - n_2 \alpha_2} I_1 + 2^{- n_1 \alpha_1} I_2 + 2^{- n_2 \alpha_2 } I_3
+  I_4.
 \end{align*}
 By Lemmas \ref{lemma1}, \ref{lemma4}    using the mixed moduli  of smoothness property, we obtain
\begin{align*}
 I_4 &\ll \Bigg( \sum\limits_{\nu_2 =n_2}^\infty \Big(
\sum\limits_{\nu_1 =n_1}^\infty \Big( 2^{\nu_1  \big( \frac{1}{p_1} - \frac{1}{q_1} \big)} 2^{ \nu_2 \big( \frac{1}{p_2} - \frac{1}{q_2} \big)  } Y_{2^{\nu_1} -1, 2^{\nu_2 }-1}(f)_{p_1 p_2} \Big)^{q_1^*}  \Big)^{\frac{q_2^*}{q_1^*}} \Bigg)^{\frac{1}{q_2^*}} \\
&  \ll \Big( \int\limits_{0}^{\delta_2} \Big( \int\limits_{0}^{\delta_1} \Big( t_1^{ - \frac{1}{p_1} + \frac{1}{q_1}} t_2^{ - \frac{1}{p_2} + \frac{1}{q_2} } \omega_{ \alpha_1  + \frac{1}{p_1} - \frac{1}{q_1}, \alpha_2  + \frac{1}{p_2} - \frac{1}{q_2}} (f, t_1, t_2)_{p_1 p_2} \Big)^{q_1^*} \frac{dt_1}{t_1} \Big)^{\frac{q_2^*}{q_1^*}} \frac{dt_2}{t_2} \Big)^{\frac{1}{q_2^*}}.
 \end{align*}
  It follows from Lemma \ref{lemma6} (b)  that
$$ I_2 \ll   \Big\| \sum\limits_{\nu_2= n_2}^\infty V_{ 2^{n_1} \infty}^{( \alpha_1, 0)} ( V_{ \infty 2^{\nu_2}}(f) -V_{\infty 2^{\nu_2 +1}}(f) ) \Big\|_{q_1q_2}. $$
Thus, by Lemma \ref{lemma5} (a), we deduce that
$$ I_2 \ll \Big( \sum\limits_{\nu_2= n_2}^\infty \Big( 2^{\nu_2 ( \frac{1}{p_2} - \frac{1}{q_2} )}  \|  V_{ 2^{n_1} \infty}^{( \alpha_1, 0)} ( V_{ \infty 2^{\nu_2}}(f) -V_{\infty 2^{\nu_2 +1}}(f) ) \|_{q_1p_2} \Big)^{q_2^*} \Big)^{\frac{1}{q_2^*}}. $$
Applying Lemma \ref{lemma12} in the case $ 1 < p_1 < q_1 < \infty  $ or Lemma \ref{lemma8} (a) in the case $ 1 = p_1 < q_1 = \infty $, we have
$$ I_2 \ll \Big( \sum\limits_{\nu_2= n_2}^\infty \Big( 2^{\nu_2 ( \frac{1}{p_2} - \frac{1}{q_2} )}  \|  V_{ 2^{n_1} \infty}^{( \alpha_1 + \frac{1}{p_1} - \frac{1}{q_1}, 0)} ( V_{ \infty 2^{\nu_2}}(f) -V_{\infty 2^{\nu_2 +1}}(f) ) \|_{p_1p_2} \Big)^{q_2^*} \Big)^{\frac{1}{q_2^*}}. $$
 Using  Lemma \ref{lemma1}, we derive the following
 \begin{align*}
&  2^{-n_1 \alpha_1} I_2 \ll \\
 &\ll  2^{n_1 (\frac{1}{p_1} - \frac{1}{q_1})} \Big( \sum\limits_{\nu_2= n_2}^\infty \Big( 2^{\nu_2 ( \frac{1}{p_2} - \frac{1}{q_2} )}\times\\
  & \;\;\;\;\;\;\;\;\;\;\;\;\;\;   \times \frac{1}{2^{n_1 (\alpha_1 + \frac{1}{p_1} - \frac{1}{q_1})}}  \|  V_{ 2^{n_1} \infty}^{( \alpha_1 + \frac{1}{p_1} - \frac{1}{q_1}, 0)} ( V_{ \infty 2^{\nu_2}}(f) -V_{\infty 2^{\nu_2 +1}}(f) ) \|_{p_1p_2} \Big)^{q_2^*} \Big)^{\frac{1}{q_2^*}} \\
  & \ll 2^{n_1 (\frac{1}{p_1} - \frac{1}{q_1})} \Big( \sum\limits_{\nu_2= n_2}^\infty \Big( 2^{\nu_2 ( \frac{1}{p_2} - \frac{1}{q_2} )} \omega_{ \alpha_1  + \frac{1}{p_1} - \frac{1}{q_1}, \alpha_2  + \frac{1}{p_2} - \frac{1}{q_2}} (f, \frac{1}{2^{n_1}}, \frac{1}{2^{\nu_2}})_{p_1 p_2} \Big)^{q_2^*} \Big)^{\frac{1}{q_2^*}} \\
    &\ll \Big( \int\limits_{0}^{\delta_2} \Big( \int\limits_{0}^{\delta_1} \Big( t_1^{ - \frac{1}{p_1} + \frac{1}{q_1}} t_2^{ - \frac{1}{p_2} + \frac{1}{q_2} } \omega_{ \alpha_1  + \frac{1}{p_1} - \frac{1}{q_1}, \alpha_2  + \frac{1}{p_2} - \frac{1}{q_2}} (f, t_1, t_2)_{p_1 p_2} \Big)^{q_1^*} \frac{dt_1}{t_1} \Big)^{\frac{q_2^*}{q_1^*}} \frac{dt_2}{t_2} \Big)^{\frac{1}{q_2^*}}.
    \end{align*}

 Now, we begin to estimate the $I_3 =  \| V_{ \infty, 2^{n_2}}^{(0, \alpha_2)} (f- V_{2^{n_1}, \infty }(f) )  \|_{q_1 q_2}.$
Applying Lemma \ref{lemma6} (a), we have
$$ I_3 \ll   \Big\| \sum\limits_{\nu_1= n_1}^\infty V_{\infty 2^{n_2}}^{(0, \alpha_2)} ( V_{2^{\nu_1} \infty}(f) - V_{2^{\nu_1 +1}\infty}(f) ) \Big\|_{q_1q_2}. $$
   Lemma \ref{lemma5} (b) yields that
$$ I_3 \ll \Big( \sum\limits_{\nu_1= n_1}^\infty \Big( 2^{\nu_1 ( \frac{1}{p_1} - \frac{1}{q_1} )}  \|  V_{\infty 2^{n_2}}^{(0, \alpha_2)} ( V_{2^{\nu_1} \infty}(f) - V_{2^{\nu_1 +1}\infty}(f) ) \|_{p_1q_2} \Big)^{q_1^*} \Big)^{\frac{1}{q_1^*}}. $$
 Applying Lemma \ref{lemma13}  in the case $1 <p_2 <q_2 <\infty $ or Lemma \ref{lemma8} (b) in the case $ 1 = p_2 < q_2 = \infty, $ we get
$$ I_3 \ll  \Big( \sum\limits_{\nu_1= n_1}^\infty \Big( 2^{\nu_1 ( \frac{1}{p_1} - \frac{1}{q_1} )}  \|  V_{\infty 2^{n_2}}^{(0, \alpha_2 + \frac{1}{p_2} - \frac{1}{q_2})} ( V_{2^{\nu_1} \infty}(f) - V_{2^{\nu_1 +1}\infty}(f) ) \|_{p_1p_2} \Big)^{q_1^*} \Big)^{\frac{1}{q_1^*}}. $$
Thus, using Lemma \ref{lemma1}, we deduce
 \begin{align*}
 & 2^{-n_2 \alpha_2} I_3 \ll \\
& \ll  2^{n_2 (\frac{1}{p_2} - \frac{1}{q_2})} \Big( \sum\limits_{\nu_1= n_1}^\infty \Big( 2^{\nu_1 ( \frac{1}{p_1} - \frac{1}{q_1} )} \times\\
& \;\;\;\;\;\;\;\;\;\; \times\frac{1}{2^{n_2 (\alpha_2 + \frac{1}{p_2} - \frac{1}{q_2})}}  \|  V_{\infty 2^{n_2}}^{(0, \alpha_2 + \frac{1}{p_2} - \frac{1}{q_2})} ( V_{2^{\nu_1} \infty}(f) - V_{2^{\nu_1 +1}\infty}(f) ) \|_{p_1p_2} \Big)^{q_1^*} \Big)^{\frac{1}{q_1^*}} \\
  &\ll 2^{n_2 (\frac{1}{p_2} - \frac{1}{q_2})} \Big( \sum\limits_{\nu_1= n_1}^\infty \Big( 2^{\nu_1 ( \frac{1}{p_1} - \frac{1}{q_1} )} \omega_{ \alpha_1  + \frac{1}{p_1} - \frac{1}{q_1}, \alpha_2  + \frac{1}{p_2} - \frac{1}{q_2}} (f, \frac{1}{2^{\nu_1}}, \frac{1}{2^{n_2}})_{p_1 p_2} \Big)^{q_1^*} \Big)^{\frac{1}{q_1^*}} \\
  &\ll \Big( \int\limits_{\frac{1}{2^{n_2+1}}}^{\frac{1}{2^{n_2}}} \Big( t_2^{ - \frac{1}{p_2} + \frac{1}{q_2} } \Big( \sum\limits_{\nu_1= n_1}^\infty \Big( 2^{\nu_1 ( \frac{1}{p_1} - \frac{1}{q_1} )} \omega_{ \alpha_1  + \frac{1}{p_1} - \frac{1}{q_1}, \alpha_2  + \frac{1}{p_2} - \frac{1}{q_2}} (f, \frac{1}{2^{\nu_1}}, \frac{1}{2^{n_2+1}})_{p_1 p_2} \Big)^{q_1^*} \Big)^{\frac{1}{q_1^*}} \Big)^{q_2^*} \frac{dt_2}{t_2} \Big)^{\frac{1}{q_2^*}} \\
   &\ll \Big( \int\limits_{0}^{\delta_2} \Big( \int\limits_{0}^{\delta_1} \Big( t_1^{ - \frac{1}{p_1} + \frac{1}{q_1}} t_2^{ - \frac{1}{p_2} + \frac{1}{q_2} } \omega_{ \alpha_1  + \frac{1}{p_1} - \frac{1}{q_1}, \alpha_2  + \frac{1}{p_2} - \frac{1}{q_2}} (f, t_1, t_2)_{p_1 p_2} \Big)^{q_1^*} \frac{dt_1}{t_1} \Big)^{\frac{q_2^*}{q_1^*}} \frac{dt_2}{t_2} \Big)^{\frac{1}{q_2^*}}.
    \end{align*}

   Let us estimate   $2^{-n_1 \alpha_1 - n_2 \alpha_2} I_1 = 2^{-n_1 \alpha_1 - n_2 \alpha_2} \| V_{2^{n_1},2^{n_2}}^{(\alpha_1, \alpha_2)}(f) \|_{q_1 q_2}.$
Applying Lemma \ref{lemma12} in the case $ 1 < p_1 < q_1 < \infty  $ or Lemma \ref{lemma8} (a) in the case  $ 1 = p_1 < q_1 = \infty  $, we obtain
$$2^{-n_1 \alpha_1 - n_2 \alpha_2} I_1 \ll 2^{-n_1 \alpha_1 - n_2 \alpha_2} \| V_{2^{n_1},2^{n_2}}^{(\alpha_1+ \frac{1}{p_1} - \frac{1}{q_1}, 0)}(f) \|_{p_1 q_2}.$$
  Using Lemma \ref{lemma13} in the case $1 < p_2 < q_2 < \infty $ or Lemma \ref{lemma8} (b) in the case $ 1 = p_2 < q_2 = \infty, $ we finally obtain
$$2^{-n_1 \alpha_1 - n_2 \alpha_2} I_1 \ll 2^{-n_1 \alpha_1 - n_2 \alpha_2} \| V_{2^{n_1},2^{n_2}}^{(\alpha_1+ \frac{1}{p_1} - \frac{1}{q_1}, \alpha_2+ \frac{1}{p_2} - \frac{1}{q_2})}(f) \|_{p_1 p_2}.$$
 By Lemma \ref{lemma1} and using the mixed moduli of smoothness property, we obtain
$$  2^{-n_1 \alpha_1 - n_2 \alpha_2} I_1 \ll
 $$
 $$
 \ll \Big( \int\limits_{0}^{\delta_2} \Big( \int\limits_{0}^{\delta_1} \Big( t_1^{ - \frac{1}{p_1} + \frac{1}{q_1}} t_2^{ - \frac{1}{p_2} + \frac{1}{q_2} } \omega_{ \alpha_1  + \frac{1}{p_1} - \frac{1}{q_1}, \alpha_2  + \frac{1}{p_2} - \frac{1}{q_2}} (f, t_1, t_2)_{p_1 p_2} \Big)^{q_1^*} \frac{dt_1}{t_1} \Big)^{\frac{q_2^*}{q_1^*}} \frac{dt_2}{t_2} \Big)^{\frac{1}{q_2^*}}. $$
 Combining the estimates for $I_1, I_2, I_3$ and $I_4$, we finally get
 $$  \omega_{\alpha_1, \alpha_2} \big( f, \delta_1,\delta_2 \big)_{q_1 q_2} \ll $$
$$ \ll \Big( \int\limits_{0}^{\delta_2} \Big(
\int\limits_{0}^{\delta_1} \Big( t_1^{ - \frac{1}{p_1} + \frac{1}{q_1}} t_2^{ - \frac{1}{p_2} + \frac{1}{q_2} } \omega_{ \alpha_1 +  \frac{1}{p_1} - \frac{1}{q_1},\alpha_2  + \frac{1}{p_2} - \frac{1}{q_2}} (f, t_1, t_2)_{p_1 p_2} \Big)^{q_1^*} \frac{dt_1}{t_1} \Big)^{\frac{q_2^*}{q_1^*}} \frac{dt_2}{t_2} \Big)^{\frac{1}{q_2^*}}. $$
Inequality \eqref{eq1} is proved.

 { \bf  Proof of Remark 1.}
Let us prove that Theorem \ref{th1} is optimal in the sense that there exists a function $f_0(x_1,x_2)$ such that for \eqref{eq1} the sign $\ll$ can be replaced by the sign $\asymp.$

  Consider the function $f_0(x_1,x_2) = \sin x_1 \sin x_2.$ For
 $\delta_i \in (0,1)$ there is a non-negative integer
a number $n_i$ such that  $\frac{1}{2^{n_1 +1}} \leq \delta_i <
\frac{1}{2^{n_i}}, i =1,2.$ Applying Lemmas \ref{lemma1} and \ref{lemma3}, for any $ r_i
> 0, i =1,2,$ and any $p_i \in [1, \infty] \; (i=1,2)$, we have
$$  \omega_{r_1, r_2} \big( f_0, \delta_1,\delta_2 \big)_{p_1,p_2} \asymp \omega_{r_1, r_2} \Big( f_0, \frac{1}{2^{n_1}}, \frac{1}{2^{n_2}} \Big)_{p_1,p_2} \asymp \frac{1}{2^{n_1 r_1 + n_2 r_2}} \asymp \delta_1^{r_1}
\delta_2^{r_2}.$$ But then
$$  \omega_{\alpha_1, \alpha_2} \big( f_0, \delta_1, \delta_2 \big)_{q_1,q_2}\asymp \delta_1^{\alpha_1} \delta_2^{\alpha_2}, $$
$$  \Big( \int\limits_{0}^{\delta_2} \Big(
\int\limits_{0}^{\delta_1} \Big( t_1^{ - \frac{1}{p_1} + \frac{1}{q_1}} t_2^{ - \frac{1}{p_2} + \frac{1}{q_2} } \omega_{ \alpha_1 +  \frac{1}{p_1} - \frac{1}{q_1},\alpha_2  + \frac{1}{p_2} - \frac{1}{q_2}} (f, t_1, t_2)_{p_1 p_2} \Big)^{q_1^*} \frac{dt_1}{t_1} \Big)^{\frac{q_2^*}{q_1^*}} \frac{dt_2}{t_2} \Big)^{\frac{1}{q_2^*}} \asymp \delta_1^{\alpha_1} \delta_2^{\alpha_2}. $$
It follows from these estimates that for the function $f_0(x_1,x_2)$ in relation \eqref{eq1} the sign $\ll$ can be replaced by the sign $\asymp.$

Another way to see the optimality of the estimate in Theorem \ref{th1} is to consider the series with monotone or general monotone coefficients, see \cite{Tikhonov4}.

{ \bf  Proof of Remark 2.}
Let us now show that, generally speaking, in inequality \eqref{eq1} the sign $\ll$ cannot be replaced by the sign $\asymp.$
  Let $ 1 < p_2 < q_2 < \infty. $ Consider the function $f_1(x_1, x_2) = \sin x_1 \cdot g(x_2),$ where $ g(x_2) = \sum\limits_{\nu =0}^\infty \frac{( \nu +1)^\beta}{2^{\nu \alpha_2}} \cos 2^\nu x_2, $  $\beta > -\frac{1}{2},  \alpha_2 > \frac{1}{p_2} - \frac{1}{q_2}.$
Since
 \begin{align*}
& \omega_{\alpha_1, \alpha_2} \big( f_1, \delta_1, \delta_2 \big)_{q_1 q_2} = \omega_{\alpha_1} \big( \sin x_1, \delta_1 \big)_{q_1}^{(1)} \cdot \omega_{ \alpha_2} \big( g, \delta_2 \big)_{ q_2}^{(1)} \asymp \delta_1^{\alpha_1} \cdot \omega_{ \alpha_2} \big( g, \delta_2 \big)_{ q_2}^{(1)}, \\
&\omega_{\alpha_1  + \frac{1}{p_1} - \frac{1}{q_1}, \alpha_2  + \frac{1}{p_2} - \frac{1}{q_2} } (f_1, t_1, t_2)_{p_1 p_2} = \omega_{\alpha_1  + \frac{1}{p_1} - \frac{1}{q_1} } (\sin x_1, t_1)_{p_1 }^{(1)} \cdot \omega_{\alpha_2  + \frac{1}{p_2} - \frac{1}{q_2} } (g,  t_2)_{ p_2}^{(1)} \\
 &\asymp t_1^{\alpha_1  + \frac{1}{p_1} - \frac{1}{q_1} }  \cdot \omega_{\alpha_2  + \frac{1}{p_2} - \frac{1}{q_2} } (g,  t_2)_{ p_2}^{(1)},
 \end{align*}
then, applying Lemma \ref{lemma14}, item (a) and substituting into inequality \eqref{eq1}, we obtain
\begin{align*}
& \omega_{\alpha_1, \alpha_2} \big( f_1, \delta_1, \delta_2 \big)_{q_1 q_2} \asymp \delta_1^{\alpha_1} \cdot \delta_2^{\alpha_2} \bigg(\log_2 \frac{2}{\delta_2} \bigg)^{\beta + \frac{1}{2}}, \\
&  \Big( \int\limits_{0}^{\delta_2} \Big(
\int\limits_{0}^{\delta_1} \Big( t_1^{ - \frac{1}{p_1} + \frac{1}{q_1}} t_2^{ - \frac{1}{p_2} + \frac{1}{q_2} } \omega_{ \alpha_1 +  \frac{1}{p_1} - \frac{1}{q_1},\alpha_2  + \frac{1}{p_2} - \frac{1}{q_2}} (f_1, t_1, t_2)_{p_1 p_2} \Big)^{q_1^*} \frac{dt_1}{t_1} \Big)^{\frac{q_2^*}{q_1^*}} \frac{dt_2}{t_2} \Big)^{\frac{1}{q_2^*}} \\
& \;\;\;\;\;\;\;\;\;\;\;\;\;\;\;\;\;
\asymp \delta_1^{\alpha_1} \cdot \delta_2^{ \alpha_2 - \frac{1}{p_2} + \frac{1}{q_2}
} \bigg( \log_2 \frac{2}{\delta_2} \bigg)^{\beta }.
\end{align*}
Thus, for the function $f_1(x_1,x_2)$, the right and left hand sides of relation \eqref{eq1} have different orders as functions of the variable $\delta_2 $ for a fixed $\delta_1$, which means that in relation \eqref{eq1} it is not possible replacing  the sign $\ll $ with the sign $\asymp.$

It can be verified similarly that in relation \eqref{eq1}, generally speaking, the sign $\ll$ cannot be replaced by the sign $\asymp$ in the case $1 < p_ 1 < q_ 1 < \infty. $

Let $1 = p_2 < q_2 = \infty.$ Consider the function $f_2(x_1, x_2) = \sin x_1 \cdot g(x_2),$ where $ g(x_2) = \sum\limits_{\nu =0}^\infty \frac{( \nu +1)^\beta}{2^{\nu \alpha_2}} \cos  \big( 2^\nu x_2 - \frac{ \pi \alpha_2}{2} \big), $   $\beta > -1, \alpha_2 > 1.$ Since
\begin{align*}
&\omega_{\alpha_1, \alpha_2} \big( f_2, \delta_1, \delta_2 \big)_{q_1 q_2} = \omega_{\alpha_1} \big( \sin x_1, \delta_1 \big)_{q}^{(1)} \cdot \omega_{ \alpha_2} \big( g_6, \delta_2 \big)_{ q_2}^{(1)} \asymp \delta_1^{\alpha_1} \cdot \omega_{ \alpha_2} \big( g, \delta_2 \big)_{ q_2}^{(1)}, \\
&\omega_{\alpha_1  + \frac{1}{p_1} - \frac{1}{q_1}, \alpha_2  + \frac{1}{p_2} - \frac{1}{q_2} } (f_2, t_1, t_2)_{p_1 p_2} = \omega_{\alpha_1  + \frac{1}{p_1} - \frac{1}{q_1} } (\sin x_1, t_1)_{p_1 }^{(1)} \cdot \omega_{\alpha_2  + \frac{1}{p_2} - \frac{1}{q_2} } (g,  t_2)_{ p_2}^{(1)} \asymp \\
 &\asymp t_1^{\alpha_1  + \frac{1}{p_1} - \frac{1}{q_1} }  \cdot \omega_{\alpha_2  + \frac{1}{p_2} - \frac{1}{q_2} } (g,  t_2)_{ p_2}^{(1)},
  \end{align*}
then, applying Lemma \ref{lemma14} item (b) and substituting into inequality \eqref{eq1}, we obtain

\begin{align*}
&\omega_{\alpha_1, \alpha_2} \big( f_2, \delta_1, \delta_2 \big)_{q_1 q_2} \asymp \delta_1^{\alpha_1} \cdot \delta_2^{\alpha_2} \bigg(\log_2 \frac{2}{\delta_2} \bigg)^{\beta + 1}, \\
& \Big( \int\limits_{0}^{\delta_2} \Big(
\int\limits_{0}^{\delta_1} \Big( t_1^{ - \frac{1}{p_1} + \frac{1}{q_1}} t_2^{ - \frac{1}{p_2} + \frac{1}{q_2} } \omega_{ \alpha_1 +  \frac{1}{p_1} - \frac{1}{q_1},\alpha_2  + \frac{1}{p_2} - \frac{1}{q_2}} (f_2, t_1, t_2)_{p_1 p_2} \Big)^{q_1^*} \frac{dt_1}{t_1} \Big)^{\frac{q_2^*}{q_1^*}} \frac{dt_2}{t_2} \Big)^{\frac{1}{q_2^*}} \asymp\\
&\;\;\;\;\;\;\;\;\;\;\;\;\;\;\;\;\;\;\;  \asymp \delta_1^{\alpha_1} \cdot \delta_2^{ \alpha_2 - 1} \bigg( \log_2 \frac{2}{\delta_2} \bigg)^{\beta + 1 }.
 \end{align*}
Therefore, for the function $f_2(x_1x_2)$, the right and left parts of relation \eqref{eq1} have different orders as functions of the variable $\delta_2 $ for a fixed $\delta_1$, which means that in relation \eqref{eq1}, it is impossible to replace the sign $\ll $ with the sign $\asymp.$

In a similar way, it is verified that, in relation to \eqref{eq1},  the sign $\ll$ cannot be replaced by the sign $\asymp$ in the case $ 1 = p_1 < q_1 = \infty. $

Thus the proof of Remark 2 is complete.

 \vskip 0.2cm

 { \bf  Proof of Theorem \ref{th2}. }
Using Theorem 1, we  have
  \begin{align*}
  &\omega_{\alpha_1, \alpha_2} \big( f^{(\rho_1, \rho_2)}, \delta_1,\delta_2 \big)_{q_1 q_2} \ll \\
 &\ll \Big( \int\limits_{0}^{\delta_2} \Big(\int\limits_{0}^{\delta_1} \Big( t_1^{ - \frac{1}{p_1} + \frac{1}{q_1}} t_2^{ - \frac{1}{p_2} + \frac{1}{q_2} } \omega_{ \alpha_1 +  \frac{1}{p_1} - \frac{1}{q_1},\alpha_2  + \frac{1}{p_2} - \frac{1}{q_2}} (f^{(\rho_1, \rho_2)}, t_1, t_2)_{p_1 p_2} \Big)^{q_1^*} \frac{dt_1}{t_1} \Big)^{\frac{q_2^*}{q_1^*}} \frac{dt_2}{t_2} \Big)^{\frac{1}{q_2^*}}.
 \end{align*}
For each $\delta_i \in (0,1)$ there is a non-negative integer $n_i$ such that $ \frac{1}{2^{n_i+1}}
\leq \delta_i < \frac{1}{2^{n_i}}, i=1,2. $ Then
 \begin{align*}
  &  \omega_{\alpha_1, \alpha_2} \big( f^{(\rho_1, \rho_2)}, \frac{1}{2^{n_1}}, \frac{1}{2^{n_2}} \big)_{q_1 q_2} \\
 &\ll \Big( \int\limits_{0}^{\frac{1}{2^{n_2}}} \Big(\int\limits_{0}^{\frac{1}{2^{n_1}}} \Big( t_1^{ - \frac{1}{p_1} + \frac{1}{q_1}} t_2^{ - \frac{1}{p_2} + \frac{1}{q_2} } \omega_{ \alpha_1 +  \frac{1}{p_1} - \frac{1}{q_1},\alpha_2  + \frac{1}{p_2} - \frac{1}{q_2}} (f^{(\rho_1, \rho_2)}, t_1, t_2)_{p_1 p_2} \Big)^{q_1^*} \frac{dt_1}{t_1} \Big)^{\frac{q_2^*}{q_1^*}} \frac{dt_2}{t_2} \Big)^{\frac{1}{q_2^*}} \\
&  \ll \Big( \sum\limits_{\nu_2=n_2}^\infty \Big( \sum\limits_{\nu_1=n_1}^\infty \Big( 2^{\nu_1 (\frac{1}{p_1} - \frac{1}{q_1})} 2^{\nu_2 (\frac{1}{p_2} - \frac{1}{q_2})} \omega_{ \alpha_1 +  \frac{1}{p_1} - \frac{1}{q_1},\alpha_2  + \frac{1}{p_2} - \frac{1}{q_2}} (f^{(\rho_1, \rho_2)}, \frac{1}{2^{\nu_1}}, \frac{1}{2^{\nu_2}})_{p_1 p_2} \Big)^{q_1^*} \Big)^{\frac{q_2^*}{q_1^*}}  \Big)^{\frac{1}{q_2^*}}. \end{align*}

a). Let  $ \rho_1>0, \rho_2=0. $ Using   Lemma \ref{lemma3}, we deduce
\begin{align*}
 &\omega_{\alpha_1, \alpha_2} \big( f^{(\rho_1, \rho_2)}, \frac{1}{2^{n_1}}, \frac{1}{2^{n_2}} \big)_{q_1 q_2} \ll\\
 & \Big( \sum\limits_{\nu_2=n_2}^\infty \Big\{ \sum\limits_{\nu_1=n_1}^\infty 2^{\nu_1 q_1^*(\frac{1}{p_1} - \frac{1}{q_1})} \times\\
 &\;\;\;\;\;\;\;  \times\Big(   \sum\limits_{m_1=\nu_1}^\infty 2^{m_1 \rho_1} \omega_{ \alpha_1 + \rho_1 +  \frac{1}{p_1} - \frac{1}{q_1},\alpha_2  + \frac{1}{p_2} - \frac{1}{q_2}} (f, \frac{1}{2^{m_1}}, \frac{1}{2^{\nu_2}})_{p_1 p_2} \Big)^{q_1^*} \Big\}^{\frac{q_2^*}{q_1^*}} 2^{\nu_2 q_2^*(\frac{1}{p_2} - \frac{1}{q_2})}  \Big)^{\frac{1}{q_2^*}}.
 \end{align*}
 Applying Lemma \ref{lemma2} to the inner sum, we obtain
\begin{align*}
 &  \omega_{\alpha_1, \alpha_2} \big( f^{(\rho_1, \rho_2)}, \frac{1}{2^{n_1}}, \frac{1}{2^{n_2}} \big)_{q_1 q_2}\ll \\
 & \Big( \sum\limits_{\nu_2=n_2}^\infty \Big\{ \sum\limits_{\nu_1=n_1}^\infty 2^{\nu_1 ( \rho_1 + \frac{1}{p_1} - \frac{1}{q_1})}     \omega_{ \alpha_1 + \rho_1 +  \frac{1}{p_1} - \frac{1}{q_1},\alpha_2  + \frac{1}{p_2} - \frac{1}{q_2}} (f, \frac{1}{2^{\nu_1}}, \frac{1}{2^{\nu_2}})_{p_1 p_2} \Big)^{q_1^*} \Big\}^{\frac{q_2^*}{q_1^*}} 2^{\nu_2 q_2^*(\frac{1}{p_2} - \frac{1}{q_2})}  \Big)^{\frac{1}{q_2^*}} \ll\\
 & \Big( \int\limits_{0}^{\delta_2} \Big(\int\limits_{0}^{\delta_1} \Big( t_1^{ - \rho_1 - \frac{1}{p_1} + \frac{1}{q_1}} t_2^{ - \rho_2 - \frac{1}{p_2} + \frac{1}{q_2} } \omega_{ \alpha_1 + \rho_1+  \frac{1}{p_1} - \frac{1}{q_1},\alpha_2 + \rho_2  + \frac{1}{p_2} - \frac{1}{q_2}} (f, t_1, t_2)_{p_1 p_2} \Big)^{q_1^*} \frac{dt_1}{t_1} \Big)^{\frac{q_2^*}{q_1^*}} \frac{dt_2}{t_2} \Big)^{\frac{1}{q_2^*}}.
 \end{align*}

b). Let $ \rho_1=0, \rho_2>0. $
By Lemma \ref{lemma3}, we have
\begin{align*}
 &\omega_{\alpha_1, \alpha_2} \big( f^{(\rho_1, \rho_2)}, \frac{1}{2^{n_1}}, \frac{1}{2^{n_2}} \big)_{q_1 q_2} \\
 &\ll \Big( \sum\limits_{\nu_2=n_2}^\infty 2^{\nu_2 q_2^*(\frac{1}{p_2} - \frac{1}{q_2})} \Big\{ \Big( \sum\limits_{\nu_1=n_1}^\infty 2^{\nu_1 q_1^*(\frac{1}{p_1} - \frac{1}{q_1})} \\
& \;\;\;\;\;\;\;\;\;\;
  \Big(   \sum\limits_{m_2=\nu_2}^\infty 2^{m_2 \rho_2} \omega_{ \alpha_1 +   \frac{1}{p_1} - \frac{1}{q_1},\alpha_2 + \rho_2  + \frac{1}{p_2} - \frac{1}{q_2}} (f, \frac{1}{2^{\nu_1}}, \frac{1}{2^{m_2}})_{p_1 p_2} \Big)^{q_1^*} \Big)^{\frac{1}{q_1^*}} \Big\}^{q_2^*}   \Big)^{\frac{1}{q_2^*}}.
 \end{align*}
Applying the Minkowski inequality to the inner sum, we get
\begin{align*}
 &  \omega_{\alpha_1, \alpha_2} \big( f^{(\rho_1, \rho_2)}, \frac{1}{2^{n_1}}, \frac{1}{2^{n_2}} \big)_{q_1 q_2}
\\
 &\ll \Big( \sum\limits_{\nu_2=n_2}^\infty 2^{\nu_2 q_2^*(\frac{1}{p_2} - \frac{1}{q_2})} \Big\{ \sum\limits_{m_2=\nu_2}^\infty  2^{m_2 \rho_2}\\
   & \;\;\;\;\;\;\;\;\;\;\; \Big( \sum\limits_{\nu_1=n_1}^\infty \Big( 2^{\nu_1 (\frac{1}{p_1} - \frac{1}{q_1})}  \omega_{ \alpha_1 +   \frac{1}{p_1} - \frac{1}{q_1},\alpha_2 + \rho_2  + \frac{1}{p_2} - \frac{1}{q_2}} (f, \frac{1}{2^{\nu_1}}, \frac{1}{2^{m_2}})_{p_1 p_2} \Big)^{q_1^*} \Big)^{\frac{1}{q_1^*}} \Big\}^{q_2^*}   \Big)^{\frac{1}{q_2^*}}.
 \end{align*}
 Furthermore, it follows from Lemma \ref{lemma2}  that
\begin{align*}
 &\omega_{\alpha_1, \alpha_2} \big( f^{(\rho_1, \rho_2)}, \frac{1}{2^{n_1}}, \frac{1}{2^{n_2}} \big)_{q_1 q_2}
\\
& \ll \Big( \sum\limits_{\nu_2=n_2}^\infty 2^{\nu_2 q_2^*( \rho_2 + \frac{1}{p_2} - \frac{1}{q_2})} \Big\{ \Big( \sum\limits_{\nu_1=n_1}^\infty \Big( 2^{\nu_1 (\frac{1}{p_1} - \frac{1}{q_1})}  \times \\
& \;\;\;\;\;\;\;\;\;\; \;\;\;\;\;\;\;\;\;\; \times \omega_{ \alpha_1 +   \frac{1}{p_1} - \frac{1}{q_1},\alpha_2 + \rho_2  + \frac{1}{p_2} - \frac{1}{q_2}} (f, \frac{1}{2^{\nu_1}}, \frac{1}{2^{\nu_2}})_{p_1 p_2} \Big)^{q_1^*} \Big)^{\frac{1}{q_1^*}} \Big\}^{q_2^*}   \Big)^{\frac{1}{q_2^*}}\ll \\
 & \Big( \int\limits_{0}^{\delta_2} \Big(\int\limits_{0}^{\delta_1} \Big( t_1^{ - \rho_1 - \frac{1}{p_1} + \frac{1}{q_1}} t_2^{ - \rho_2 - \frac{1}{p_2} + \frac{1}{q_2} } \omega_{ \alpha_1 + \rho_1+  \frac{1}{p_1} - \frac{1}{q_1},\alpha_2 + \rho_2  + \frac{1}{p_2} - \frac{1}{q_2}} (f, t_1, t_2)_{p_1 p_2} \Big)^{q_1^*} \frac{dt_1}{t_1} \Big)^{\frac{q_2^*}{q_1^*}} \frac{dt_2}{t_2} \Big)^{\frac{1}{q_2^*}}.
 \end{align*}

c). Let $ \rho_1>0, \rho_2>0. $ Using Lemma \ref{lemma3}, we obtain
\begin{align*}
 &  \omega_{\alpha_1, \alpha_2} \big( f^{(\rho_1, \rho_2)}, \frac{1}{2^{n_1}}, \frac{1}{2^{n_2}} \big)_{q_1 q_2} \ll  \Big( \sum\limits_{\nu_2=n_2}^\infty 2^{\nu_2 q_2^*(\frac{1}{p_2} - \frac{1}{q_2})} \Big\{ \Big( \sum\limits_{\nu_1=n_1}^\infty 2^{\nu_1 q_1^*(\frac{1}{p_1} - \frac{1}{q_1})} \times
\\
 & \times  \Big\langle   \sum\limits_{m_1=\nu_1}^\infty 2^{m_1\rho_1}  \sum\limits_{m_2=\nu_2}^\infty 2^{m_2 \rho_2} \omega_{ \alpha_1 + \rho_1 + \frac{1}{p_1} - \frac{1}{q_1},\alpha_2 + \rho_2  + \frac{1}{p_2} - \frac{1}{q_2}} (f, \frac{1}{2^{m_1}}, \frac{1}{2^{m_2}})_{p_1 p_2} \Big\rangle^{q_1^*} \Big)^{\frac{1}{q_1^*}} \Big\}^{q_2^*}   \Big)^{\frac{1}{q_2^*}}.
 \end{align*}
 Applying Lemma \ref{lemma2} to  the sum in the "angle" brackets,  we get
\begin{align*}
  &\omega_{\alpha_1, \alpha_2} \big( f^{(\rho_1, \rho_2)}, \frac{1}{2^{n_1}}, \frac{1}{2^{n_2}} \big)_{q_1 q_2} \ll  \Big( \sum\limits_{\nu_2=n_2}^\infty 2^{\nu_2 q_2^*(\frac{1}{p_2} - \frac{1}{q_2})} \Big\{ \Big( \sum\limits_{\nu_1=n_1}^\infty 2^{\nu_1 q_1^*( \rho_1 + \frac{1}{p_1} - \frac{1}{q_1})} \times \\
  &\times  \Big\langle     \sum\limits_{m_2=\nu_2}^\infty 2^{m_2 \rho_2} \omega_{ \alpha_1 + \rho_1 + \frac{1}{p_1} - \frac{1}{q_1},\alpha_2 + \rho_2  + \frac{1}{p_2} - \frac{1}{q_2}} (f, \frac{1}{2^{\nu_1}}, \frac{1}{2^{m_2}})_{p_1 p_2}  \Big\rangle^{q_1^*} \Big)^{\frac{1}{q_1^*}} \Big\}^{q_2^*}   \Big)^{\frac{1}{q_2^*}}.
  \end{align*}
  By application of Minkowski's inequality, we obtain
\begin{align*}
  &\omega_{\alpha_1, \alpha_2} \big( f^{(\rho_1, \rho_2)}, \frac{1}{2^{n_1}}, \frac{1}{2^{n_2}} \big)_{q_1 q_2} \ll  \Big( \sum\limits_{\nu_2=n_2}^\infty 2^{\nu_2 q_2^*(\frac{1}{p_2} - \frac{1}{q_2})} \Big\{ \sum\limits_{m_2=\nu_2}^\infty  2^{m_2 \rho_2} \times
\\
& \times  \Big( \sum\limits_{\nu_1=n_1}^\infty \Big\langle 2^{\nu_1 ( \rho_1 + \frac{1}{p_1} - \frac{1}{q_1})}       \omega_{ \alpha_1 + \rho_1 + \frac{1}{p_1} - \frac{1}{q_1},\alpha_2 + \rho_2  + \frac{1}{p_2} - \frac{1}{q_2}} (f, \frac{1}{2^{\nu_1}}, \frac{1}{2^{m_2}})_{p_1 p_2}  \Big\rangle^{q_1^*} \Big)^{\frac{1}{q_1^*}} \Big\}^{q_2^*}   \Big)^{\frac{1}{q_2^*}}.
\end{align*}
using Lemma \ref{lemma2} again, we deduce that
\begin{align*}
  & \omega_{\alpha_1, \alpha_2} \big( f^{(\rho_1, \rho_2)}, \frac{1}{2^{n_1}}, \frac{1}{2^{n_2}} \big)_{q_1 q_2} \ll  \Big( \sum\limits_{\nu_2=n_2}^\infty 2^{\nu_2 q_2^*( \rho_2 + \frac{1}{p_2} - \frac{1}{q_2})} \times
\\
& \times \Big\{  \Big( \sum\limits_{\nu_1=n_1}^\infty \Big\langle 2^{\nu_1 ( \rho_1 + \frac{1}{p_1} - \frac{1}{q_1})}       \omega_{ \alpha_1 + \rho_1 + \frac{1}{p_1} - \frac{1}{q_1},\alpha_2 + \rho_2  + \frac{1}{p_2} - \frac{1}{q_2}} (f, \frac{1}{2^{\nu_1}}, \frac{1}{2^{\nu_2}})_{p_1 p_2}  \Big\rangle^{q_1^*} \Big)^{\frac{1}{q_1^*}} \Big\}^{q_2^*}   \Big)^{\frac{1}{q_2^*}}\ll \\
& \Big( \int\limits_{0}^{\delta_2} \Big(\int\limits_{0}^{\delta_1} \Big( t_1^{ - \rho_1 - \frac{1}{p_1} + \frac{1}{q_1}} t_2^{ - \rho_2 - \frac{1}{p_2} + \frac{1}{q_2} } \omega_{ \alpha_1 + \rho_1+  \frac{1}{p_1} - \frac{1}{q_1},\alpha_2 + \rho_2  + \frac{1}{p_2} - \frac{1}{q_2}} (f, t_1, t_2)_{p_1 p_2} \Big)^{q_1^*} \frac{dt_1}{t_1} \Big)^{\frac{q_2^*}{q_1^*}} \frac{dt_2}{t_2} \Big)^{\frac{1}{q_2^*}}.
\end{align*}
Theorem \ref{th2} is completely proved.

\


\begin{thebibliography}{99}

\bibitem {Besov1}
O.V. Besov, V.P. Ilyin, S.M. Nikolsky, \textit{  Integral
representation of functions and the embedding theorem}, J. Wiley and Sons, New York, 1978; translated from Russian: Nauka, Moscow, 1975.


\bibitem {DitTik}
Z. Ditzian, S. Tikhonov, \textit{ Ul'yanov and Nikol'skii-type inequalities}, Journal of Approx. Theory 133(1) (2005), 100--133.

\bibitem {DomTik}
 O. Domingues, S. Tikhonov, \textit{ Embedding of smooth function spaces extrapolations, and  related  inequalities,} ArXiv: 1909.12818v2 [math.FA] (2019), 1--71.



\bibitem {GlazTik}
 P. Glazyrina,  S. Tikhonov, \textit{ Jacobi weights, fractional integration,
and sharp Ulyanov inequalities},  J. Approx. Theory  195  (2015), 122--140.

\bibitem {GogOp}
A. Gogatishvili,  B. Opic, S. Tikhonov, W. Trebels,  \textit{ Ulyanov-type  inequalities between Lorentz-Zygmund spaces}, J. Fourier Anal. Appl. 20(5) (2014), 1020--1049.


\bibitem {Goldm} M.L. Gol'dman,  \textit{  Embedding of constructive and structural Lipschitz spaces in symmetric spaces}, Tr. Mat. Inst. Steklova 173 (1986), 90--112 (in Russian); translated in: Proc Steklov Inst. Math. 173(4) (1987), 93--118.

\bibitem{Ainur}
A.A. Jumabayeva, \textit{ Sharp Ul'yanov inequalities for generalized Liouville-Weyl derivatives}, Anal. Math. 43(2) (2017),  279–302.


\bibitem{JumSim}
 A.A. Jumabayeva,  B.V. Simonov, \textit{   Transformation of Fourier series by means of general monotone sequences},   Mat. Zametki 107(5) (2020), 674--692.


\bibitem {Kolyada}
 V.I. Kolyada,  \textit{ On the relations between the moduli of continuity in different metrics,} Trudy Mat. Institute of the Academy of Sciences of the USSR. 181 (1988), 117--136.

 \bibitem {KolTik1}
  Yu. Kolomoitsev, S. Tikhonov,  \textit{ Properties of moduli of smoothness in
$L_p(\mathbb{R}^d)$,} J. Approx. Theory 257 (2020), Article  105423.
Arxiv: 1907.12788, 1--29.



 \bibitem {KolTik3}  Yu. Kolomoitsev, S. Tikhonov, \textit{ Hardy-Littlewood and Ulyanov inequalities,} Mem Amer. Soc. 271(1325) (2021),  Arxiv: 1711.08163.




\bibitem {Potapov2} M.K. Potapov, B.V. Simonov, S. Yu. Tikhonov, \textit{ Mixed moduli of
smoothness in $L_p,  1 < p < \infty:$  a survey},  Surveys in
Approximation Theory, 8 (2013), 1--57.


\bibitem {Potapov3} M. K. Potapov, B. V. Simonov, S. Yu. Tikhonov, \textit{  Relations between mixed moduli of smoothness and embedding theorems for Nikol'skii classes}, Proc. Steklov Inst. Math., 269 (2010), 197–207


\bibitem {Potapov4} M. K. Potapov, B. V. Simonov, S. Yu. Tikhonov, \textit{ Analogues of Ulyanov Inequalities for Mixed Moduli of Smoothness}, Methods of Fourier Analysis and Approximation Theory, Applied and Numerical Harmonic Analysis, Springer International Publishing, Switzerland, 2016, 161--179.

\bibitem {Potapov5}  M. K. Potapov, B. V. Simonov, \textit{ Refinement of the relations between mixed smoothness moduli in $L_1$  and $L_q$ metrics}, Siberian Math. J. 62(4) (2021), 661--677.

\bibitem {Potapov6} M. K. Potapov, B. V. Simonov, \textit{ Refinement of Relations between Mixed Moduli of Smoothness in the Metrics of $L_p$ and $L_\infty$},  Math. Notes 110(3) (2021), 347--362.



\bibitem {Potapov7} M. K. Potapov, B. V. Simonov, \textit{ Properties of mixed moduli of smoothness of positive order in a mixed metric}, Vestn. Mosk. un-ta. Ser.  Matem. mekh. 6 (2014),  31--40.


\bibitem {Potapov8} M. K. Potapov, \textit{ Embedding theorems in a mixed metric}, Tr. MIAN SSSR 156 (1980), 143--156.


\bibitem {PS8}
 M.K. Potapov, B.V. Simonov,  \textit{ Inequalities of different metrics for trigonometric polynomials}, Izvestiya vuzov. Math. 1 (2019),  49--62.


\bibitem {PS6}
M.K. Potapov, B.V. Simonov, \textit{ Estimates of partial moduli of smoothness in metrics of $L_{p_1 \infty}$ and $L_{\infty p_2}$ through partial moduli of smoothness in metrics of $L_{p_1 p_2},$} Moscow Univ. Math. Bull.  75 (2020), 3--17.


 \bibitem {PST4}
    M. K. Potapov, B. Simonov, S. Tikhonov, \textit{  Fractional moduli of smoothness}, Maks--Press, Moscow, 2016.




\bibitem {PS5}
  M. K. Potapov, B.V. Simonov, \textit{ Strengthened Ul'yanov's inequalities for partial moduli of smoothness for functions from spaces with various metrics,} Vestnik Moskov. Univ. Ser.  Mat. Mekh. 3 (2019),  26--38.


\bibitem {ST3}
 B. Simonov, S. Tikhonov, \textit{ Sharp Ul'yanov-type inequalities
using fractional smoothness}, J. Approx. Theory, 162(9) (2010),   1654--1684.

\bibitem {Tikhonov1}
 S. Tikhonov,  \textit{ Weak type inequalities for moduli of smoothness:
the case of limit value parameters,} J. Fourier Anal.
Appl. 16(4) (2010),  590--608.



  \bibitem {Tikhonov3}
   S. Tikhonov,  \textit{ Trigonometric series of Nikol'skii classes},  Acta Math. Hungar. 114(1-2)  (2007), 61--78.

 \bibitem {Tikhonov4}
 S. Tikhonov,  \textit{ Trigonometric series with general monotone coefficients},
Jour. Math. Anal. Appl. 326(1) (2007), 721--735.


\bibitem {TikTre}
S. Tikhonov, W. Trebels, \textit{ Ulyanov inequalities and generalized Liouville derivatives}, Proc. Roy. Soc. Edinburgh Sect. A  141(1) (2011),  205--224.



\bibitem {Trebels}  W. Trebels, \textit{ Inequalities for moduli of smoothness versus
embeddings of function spaces,} Arch. Math. 94 (2010), 155--164.


  \bibitem {Ul1}
 P.L. Ul'yanov,  \textit{ The imbedding of certain function classes $H_p^{\omega}$},  Izv. Akad. Nauk SSSR Ser. Mat. 32(3) (1968),  649--686.

\bibitem {Zigmund}
A. Zygmund, \textit{ Trigonometric  Series,}  Vol. I. II Third edition. Cambridge 2002.


\end{thebibliography}
                         \end{document}